\documentclass[10pt,a4paper,oneside,onecolumn]{article}

\usepackage{enumerate}

\usepackage{mathtools}
\usepackage[T1]{fontenc}
\usepackage{amsmath}
\usepackage{amsfonts}
\usepackage{amssymb}
\usepackage{epsfig}
\usepackage{makeidx}
\usepackage{graphicx}
\usepackage{caption}
\usepackage{subcaption}
\usepackage{xcolor}
\usepackage{tcolorbox}
\usepackage{float}
\usepackage{tabulary}
\usepackage{tabularx}
\usepackage{booktabs}
\usepackage{color}
\usepackage{authblk}
\usepackage{setspace}
\usepackage{bm}
\usepackage{comment}
\usepackage{enumitem}   
\usepackage{ upgreek }
\usepackage{tikz}
\usepackage{wrapfig}
\usepackage{framed}
\definecolor{mycolor}{rgb}{0.01,0.1,0.3}
\usetikzlibrary{shapes.geometric, snakes, shapes.arrows, patterns, shapes.multipart, mindmap, trees, calc, shapes, decorations.text}
\tikzset{concept/.append style={fill={none}}}

\usepackage{todonotes}
\usepackage{siunitx} % Required for alignment
\allowdisplaybreaks

\usepackage{algorithm}
\usepackage{algorithmicx}
\usepackage{algpseudocode}

\definecolor{shadecolor}{rgb}{0.98,0.9,1}

\usetikzlibrary{arrows.meta}

\providecommand{\keywords}[1]{\textbf{\textit{Keywords---}} #1}
\restylefloat{figure}
\title{Optimization }
\author[1]{Maurizio Boccia}
\author[2]{Veronica Dal Sasso}
\author[3]{Leonardo Lamorgese}
\author[4]{Carlo Mannino}
\author[5]{Paolo Ventura}
\affil[1]{University of Naples, Italy,  maurizio.boccia@unina.it}
\affil[2]{Siemens Mobility, Rome, Italy, veronica.dal-sasso@siemens.com}
\affil[3]{Siemens Mobility, Rome, Italy, leonardo.lamorgese@siemens.com}
\affil[4]{Siemens Mobility, Norway, carlo.mannino@siemens.com}
\affil[5]{Siemens Mobility, Rome, Italy, paolo.ventura@siemens.com}

\newcommand{\CAR}[1]{\todo[inline, color=orange!40]{CAR: #1}}

\newcommand{\reals}{\mbox{${\textrm I}\!{\textrm R}$}}

\newtheorem{theorem}{Theorem}
\newtheorem{problem}[theorem]{Problem}
\newtheorem{lemma}[theorem]{Lemma}

\begin{document}

\title{Optimizing train dispatching for the Union Pacific Railroad}


%\KEYWORDS{deterministic inventory theory; infinite linear programming duality; 
%  existence of optimal policies; semi-Markov decision process; cyclic schedule}
%\HISTORY{This article was reviewed} % for example

% Fill in data. If unknown, or unnecessary, outcomment the field
\maketitle
\begin{abstract}
Union Pacific (UP) is one of the largest transportation companies in the world, with over 50.000 kms of rail network covering 23 states in the United States. In 2017 Union Pacific embarked on a project that within 5 years would lead it to become the only rail operator in the world equipped with a technology capable of fully automating the real-time management and optimization of train traffic. In 2021 the main milestone of such project has been reached with the first deployment of the automated dispatching system we present here.
To attack such large and complex problem, we decomposed it into distinct but interrelated functional components, and developed optimization models and methods to handle such components. The models communicate with each other through variables and constraints, and by a careful timing of invocations.  In this paper we give an overview of the overall approach.  % Enter your abstract
\end{abstract}

\keywords{Optimization, traffic management, A-TMS, real-life application}

\section{Introduction}\label{sec:introduction}

Union Pacific owns and manages a large and complex rail network with hundreds or even thousands of trains per day, hauling goods and, to a lesser extent, carrying passengers from one side to the other of the country. Like with most rail networks, the tracks are segmented into elementary rail resources (the {\em block sections}) which, for safety reasons, can generally be occupied by at most one train at a time. Since the main objective is for trains to reach their destination in as little time as possible, it becomes clear how questions such as “which train goes first?” or “is an alternative route available?” must be answered over and over by the controllers in charge of the network, the so-called dispatchers. These decisions are very challenging to take, because dispatchers have relatively little time to take them (sometimes only a few seconds) and because each decision can, potentially, have a cascading effect on the rest of the network. Factoring in the implications of each decision before taking it is extremely hard, even for experienced and skilled dispatchers. What might appear the best decision locally can create unforeseen delays or issues down the line, even in different parts of the network. Furthermore, in the North-American context this challenge is compounded by two important factors. 
First of all, trains are often "overlength", in the sense that they are longer than the majority of rail resources. This drastically increases the risk of creating a {\em deadlock}, which, broadly speaking, is a situation in which a subset of trains are blocking mutually needed resources and no feasible move allows them to reach their destination. Deadlocks are very pernicious events for a rail operator, as the recovery process is typically complex and very costly. 
The second factor is that UP, like all North-American railroads, traditionally operates "unscheduled", that is, with no official timetable to adhere to.\footnote{Unlike networks primarily focused on passenger traffic, such as those common in Europe} This makes a dispatcher's job arguably even harder, as train composition and traffic patterns change continuously. 

In this setting, two strategic decisions were taken at UP that were to have profound consequences. The first was the complete renovation of its Traffic Management system (TMS - also known as Dispatch system). A TMS is the command-and-control software interface between the dispatchers and signalling system, which allows centralised control of train movements by means of activating switches, signals, etc. The novelty lies in the fact that the renovation was done entirely {\em in-house}, the first rail operator to ever do so (at least to the best of our knowledge). The second, crucial decision was to determine that the TMS should be equipped with planning technology to fully automate the process of managing and optimizing the network and its traffic. To achieve this, UP decided to integrate its TMS with AMP, a planning software that harnesses different optimization algorithms to produce real-time scheduling and routing plans for the entire controlled network, detecting and avoiding deadlocks, coordinating traffic between different regions, and more, all in real-time. Effectively, AMP is an "artificial dispatcher", providing decision support to human dispatchers and taking full control when launched in auto-mode. 

These two decisions would lead to the successful integration and implementation in May 2021 of this innovative technology, which we will refer to as an A-TMS (or Automated Traffic Management system). System-wide rollout was then completed by November of the same year. 
In this paper we give an overview of the variety of complex problems tackled by AMP. In particular, we give a hint of the mathematics that was developed to solve them, that ranges from graph theory, combinatorial optimization, exact and heuristic methods for integer linear programming and decomposition theory. We refer to the various papers for more detail on the approaches. 

\subsection{An overview of the real-time train traffic management problem} \label{problem_overview}
The problem of managing, or {\em dispatching}, trains in real-time has presumably existed since the inception of railway systems in the nineteenth century. Having to sequence and schedule trains on constrained resources requires a level of coordinated and (ideally) informed decision-making. As the scale and complexity of rail networks has developed over time, so, fortunately, has the technology used to govern them. Train traffic around the world is nowadays controlled by signals, interlocking and so-called Automatic Train Protection systems, thanks to which crucial tasks such as setting train routes and imposing safety distance between trains can be carried out in a safe, informed and even centralised way via the TMS.
The safety of {\em fixed block} systems is ensured by allowing at most one train at the time on certain "atomic" resources, such as station platforms or line block sections (as mentioned previously, portions of track delimited by two consecutive signals). So-called {\em moving block} systems rely instead on defining a real-time safety zone around each train, and ensuring that no other trains enter such area.\footnote{This requires exact knowledge of the position and speed of all trains and is therefore a technologically more advanced solution} 
The path of a train in the network can be seen as a sequence of rail resources, each of which is occupied by a train for a given time interval that we refer to as {\em occupation time}. This time represents the time elapsing between the entrance of the train's head and the exit of the train's tail (plus some additional time to release the resource), and depends on a number of things such as train length, speed, potential dwell and more.\footnote{The concept in moving block systems is slightly more complex so we omit it for sake of brevity}
The main point is that a rail network is heavily resource-constrained, since these occupation times must be accounted for when trains travel towards their destination. It is therefore not hard to see how delays accumulate rapidly, as the primary delay of a train due to an unforeseen circumstance may easily propagate to other trains. Indeed, the major challenge in managing a rail network is the constant need to take routing, sequencing and scheduling decisions to minimize the impact of these delays. This is why the dispatcher's role is so arduous but, at the same time, crucial.

A rail network of large size may require tens or even more than a hundred dispatchers constantly overseeing operations. The network is generally partitioned into different areas, each controlled by a regional operating center and/or by a group of dispatchers, where in turn each dispatcher is assigned a smaller portion of this region, which, depending on average traffic conditions, could be as granular as a line, a station or even more microscopic. 
%eliminabile
{This is somewhat analogous to the way air traffic operations are handled, where traffic controllers in a specific tower "take over" the supervision of a plane once it enters their airspace and accompany it end-to-end until it leaves the airspace (including potential landing, taking off and ground movements in the airport). }

\begin{wrapfigure}{}{0.5\textwidth}
%\begin{figure}
%    \centering
    \includegraphics[width=0.5\textwidth]{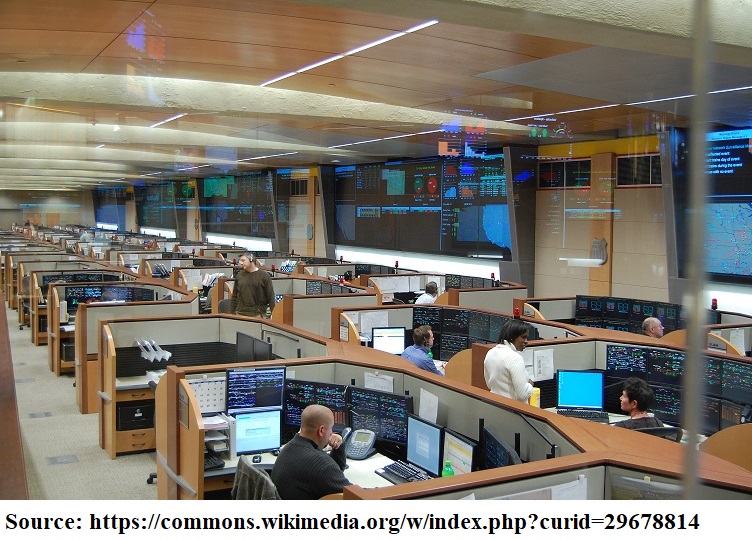}
    \caption{A dispatcher operating center}
    \label{fig:harriman}
%\end{figure}
\end{wrapfigure}

This {\em divide et impera} approach allows dispatchers to tame the complexity of the task they face, at least to some extent. On the other hand, it may lead to taking suboptimal choices from a global perspective. Indeed, each dispatcher naturally focuses primarily on dispatching her area of competence as well as possible. While there is of course some level of coordination between dispatchers (especially those that sit in the same operating centre), it is virtually impossible for any human being to fully assess the impact of each decision on the entire network, particularly when actively concentrating only on a small portion of it. An example of an operating center is given in Figure \ref{fig:harriman}, where dispatchers can be seen at work.

\subsection{A-TMS state-of-the-practice} 
The primary goal of an automated TMS is to find (typically in a few seconds) a schedule and routing - that is, {\em a plan} - for trains in a railway network, for the next $n$ hours - with $n$ decided by the rail operator. Such a plan must be conflict- and deadlock-free, in accordance with the destination(s) of the trains. Moreover, the plan should minimize some measure of the delay of the trains. As such, the addressed problem is an  optimization problem, often referred to as {\em train rescheduling} problem, or simply {\em train dispatching} problem. There is a vast body of literature, which can be classified in various ways. Since the focus of this section is to give an overview of the state of the practice, we classify here the previous contributions according to the degree of implementation in real-life - and not, for instance, by the mathematical model or solution methodology. In this classification, we start from purely "theoretical" works, which test their algorithms on (1) {\em artificial instances}; we then have methods which are tested on (2) {\em realistic or real-life instances}, (3) {\em pilot studies} which are tested on-field and (4) {\em real-life implementations}, which are currently in operation or have been at some point in time. In this section we focus only on (3) and (4), and refer to various surveys and volume articles (such as \cite{cacchiani2014overview, fang2015survey, lamorgese2018train, marcelli2020literature, wen2019train}) for the remaining works. Note that in contrast with (1) and (2), the literature on (3) and (4) is very limited and typically lacking of crucial details to allow a full understanding of the chosen approach.

\paragraph{Pilot studies.} In our definition of a pilot study, the optimization tool must have been tested by dispatchers in a real-life production environment, with data (such as the position of the trains and the status of the infrastructure) continuously updated from the field. The main difference with a real-life implementation is that the algorithm is not integrated with the TMS, and therefore the plans produced are not directly displayed to the dispatchers in the system and, obviously, cannot be carried out directly in the field. Pilot studies also tend to focus on smaller portions of a network that present some properties that make them suitable for a first tussle with real-life. 

The earliest pilot study is described in \cite{mazzarello2007traffic} and was carried out in 2004 by the Dutch infrastructure manager Pro-Rail, within the  project COMBINE 2 \cite{giannettoni2004european}.  The pilot area in the Netherlands started at Breda and at Roosendaal and ended at Willemsdorp. The network was a single line with a bifurcation. Interestingly, this early study already showed that the average punctuality of the system could be improved following the recommendations of an optimization algorithm.

An interesting case is the approach described in \cite{bettinelli2017real}, as this was later embedded in a commercial tool (ICONIS, by Alstom). Nevertheless, to the best of our knowledge, the automatic planning functions have never been deployed in a production setting, but limited to testing in pilot experiments in France.\footnote{personal communication}  

A very recent experience is the optimization train rescheduling system \cite{Brugmann,Hostettler2020}  developed in-house at  SBB, the Swiss state-held railway infrastructure and train manager. SBB developed a Branch\&Bound algorithm which is now undergoing a pilot test-campaign on field.   The tests are carried out in two regions: (a) the area around Olten where many different lines come together and part again (trains from/to Zurich, Basel, Bern, Biel, Lucerne, and other cities); 
(b)	the area between Bellinzona, Lugano and Locarno, including the new Ceneri Base Tunnel (CBT).
It handles up to 60-80 trains, for a planning horizon of 15 to 30 minutes.

\paragraph{Real-life implementations.} Finally, we describe the few developments which are in operation or have been in operation on some railways. The plans returned by the tool may be automatically transmitted to the field, or simply visualized to dispatchers which can accept or modify the tools decisions. The first known example of a real-life implementation of an optimized TMS is the control system of the Lötschberg Base Tunnel (SBB/Syntransis) control system \cite{montigelsemi} in 2006. However, from the limited descriptions available, somehow confirmed in \cite{rao2016new}, the system is more akin to a speed control system.

Next, a full-fledged automatic TMS was put in operation by Bombardier Transportation in year 2007 in one of the stations of Milan metro system \cite{mannino2009optimal}, based on an exact Branch\&Bound approach. Plans are generated and returned every 10 secs. Despite the relatively small size of the station, the system performed on average better than the dispatchers on a test-campaign carried out prior to the actual commissioning. The system was equipped with automatic route setting functions. The system was switched off when Bombardier Transportation lost the tender for the complete renewal of the train control system in Milan underground. 

Bombardier Transportation later commissioned mainline A-TMSs in 2011 and 2012 for some Italian regional lines (and later in 2016 some in Latvia), see \cite{lamorgese2018train}. All these systems implement a heuristic version of the (Benders' like) decomposition approach described in \cite{lamorgese2015exact}. Plans are generated and returned every 10 secs, and the number of controlled trains is around 100, for a planning horizon of 6 hours. Although equipped with automatic route setting functions, the plans are not sent automatically to  the field, but rather presented to dispatchers which must confirm or refuse the meet/pass solutions proposed by the system. 

In 2012 General Electric (GE) commissioned an A-TMS (referred to as {\em movement planner}) for the American class I railway Norfolk Southern. The experience and, briefly, the approach are described in \cite{bollapragada2018novel}. The approach is a combination of several optimization methods, embedded into a beam-search Branch\&Bound. From public documentation, it is hard to infer the actual degree of automation, the feasibility and quality of the delivered plans, the cycle time of the system and other crucial details. Nevertheless, the GE movement planner (now acquired by Wabtec) represents a remarkable achievement in the development and deployment of A-TMS worldwide. 

\subsection{Introducing the next sections}
In the next section we present a well-known model for train dispatching. The following section gives an overview of the AMP system, including the description of the core train dispatching functionality, the ancillary functionalities of deadlock avoidance and deadlock detection, and techniques to tackle scalability. The aim of this paper is to describe the role of the individual components in such complex system, and to provide some details about the techniques adopted to solve the individual subproblems. For more details on the mathematics developed to tackle these problems and associated computational results we will refer, when possible, to published papers.

\section{A model for the train dispatching problem} 
It is well established that train scheduling can be modelled as a job-shop scheduling problem (with routing) with blocking and no-wait constraints \cite{mascis2002job}, where trains correspond to jobs and tracks correspond to machines. As such, an instance of train scheduling can be represented by a suitable disjunctive graph \cite{balas1969machine}, which plays a central role in our approach.  

In this section we will briefly describe how this graph is built  from an instance of train dispatching. We need to first represent the possible movements of a set $T$ of trains in the railway network. To this end, the tracks of the network are subdivided into previously introduced smaller segments, called {\em block sections}. For safety reasons, each block section can accommodate at most one train at a time. 
A train runs from its origin to its destination through a sequence of block sections, called the {\em route} of the train. A route can be modeled as a graph, namely by a directed path, with arcs corresponding to block sections\footnote{This is a simplification for the sake of simplicity, as in the actual model the arcs represent {\em sequences} of block sections. Similarly, we assume that any train can be entirely contained in a block section.} and nodes corresponding to the junction between adjacent sections. There are in general multiple routes available for each train. The set of routes associated with train $i\in T$ is represented by a directed graph $G_i = (N_i, A_i)$, with a source node $o_i\in N_i$ corresponding to the origin of the train and a sink node $d_i\in N_i$ corresponding to its destination. More in general, the nodes $N_i$ correspond to all block sections the train $i$ could traverse while each arc $(u,v) \in A_i$ represents the fact that $i$ can enter $v$ from $u$. $G_i$ is the {\em routing graph} of train $i$
By construction, for any $v\in N_i$, there is at least one directed path from $o_i$ to $v$  and from $v$ to $d_i$. Note that choosing a route for train $i$ corresponds to choosing a path from $o_i$ to $d_i$ in $G_i$, that is an ordered set of arcs $P = ((v_0, v_1)\dots (v_{q-1},v_q))$ of $G_i$, with $v_0 = o_i$ and $v_q = d_i$. 

In a slightly different perspective, we can look at $G_i$ as a {\em discrete event} graph. Under this interpretation, each node represents the {\em event of entering a block section}. The arc $(u,v)\in A_i$ represents the fact that, for train $i$, $u$ and $v$ are associated with adjacent block sections and the event of entering the block section (associated with) $v$ occurs {\em after} the event of entering the block section (associated with) $u$. Therefore, we associate with arc $(u,v)\in A_i$ the length $l_{uv}\in \reals_+$ to represent the minimum running time of train $i$ to traverse the block section (associated with) $u$. Given a feasible route $P = ((v_0, v_1)\dots (v_{q-1},v_q))$ for train $i$, a {\em schedule} for $i$ is an assignment of times $t_0, \dots, t_q \in \reals_+$ to the nodes of $P$ satisfying $t_k-t_{k-1} \geq l_{k-1,k}$, for $k = 1,\dots, q$. 

Now, once a route is established for each train, we must decide how trains are sequenced on contended resources. 
To this end, we introduce a new event graph $G=(N,A)$ (see \cite{mannino2020exact}), which is the union of all routing graphs, plus some arcs and one extra node $s$. The extra node $s$  represents the event "start of the planning horizon" and we have $N = \{s\} \cup_{i\in T} N_i$. Note that, for a given block section $b$, if $u_i \in N\cap N_i$ is the event {\em train $i$ enters $b$} and $u_j\in N\cap N_j$ is the event {\em train $j$ enters $b$}, $u_i$ and $u_j$ are distinct nodes of $N$. 

The set $A$ is partitioned into three subsets: the {\em routing arcs} $A^R = \cup_{i\in T} A_i$, which is simply the union of all arcs in the routing graphs; the {\em timetable arcs} $A^Q$ and the {\em conflict arcs} $A^C$. Timetable arcs $A^Q$ are all outgoing the start node $s$ and represent lower bounds on the time the events can happen. In particular, the timetable arc $(s,u)\in A^Q$ with length $l_{su}\in \reals_+$ implies the event associated with node $u$ (for instance a train departure from a station) cannot occur earlier than $l_{su}$.  We assume $(s,o_i) \in A^Q$, for $i\in T$, with $l_{s{o_i}} \geq 0$, representing that every train starts its movement not earlier than the beginning of the planning horizon.  Conflict arcs $A^C$ model the sequencing on contended block sections. Suppose $u_i, u_j\in A$ represent train $i$ and train $j$ entering a same physical block section $b$. This situation is called a (potential) {\em  conflict}. Let $v_i$ ($v_j$) represent train $i$ (train $j$) entering the block section immediately following $b$ in its route. Then, either train $i$  traverses $b$ entirely before train $j$ enters $b$ ($i$ {\em wins} the conflict), or viceversa ($j$ wins the conflict). If train $i$ wins, then $i$ enters the next block before $j$ enters $b$, and this can be represented by the (conflict) arc $(v_i,u_j)$ with length $l_{v_i,u_j} = 0$. If $j$ wins the conflict, then $j$ runs first through $b$, represented by the (conflict) arc $(v_j,u_i)$ with length $l_{v_j,u_i} = 0$. The pair of conflict arcs $\{(v_i,u_j),(v_j,u_i)\}$ is called {\em disjunctive couple} \cite{balas1969machine}, and is associated with the conflict of trains $i$ and $j$ on block section $b$. Deciding which of the two trains wins the potential conflict is equivalent to choosing one arc in the disjunctive couple. 
Observe that
\begin{itemize}
    \item Deciding a route for each train amounts to choosing a suitable subset  of the routing arcs $A^R$, namely a set of arcs defining an origin destination path in $G^i$, for each $i\in T$.
    \item Given the route of each train, deciding which train wins each potential conflict amounts to choosing a suitable subset of conflict arcs $A^C$, namely one arc for each potential conflict (disjunctive couple) implied by the chosen routes.  
\end{itemize}

So, an instance of the train dispatching problem can be described by means of graph $G=(N,A)$ above defined, and we call it {\em dispatching graph}. As we have seen, finding a route for all trains and deciding how to sequence the trains on the contended block sections amounts to choosing a suitable subset $y\subseteq A$ of arcs. We call $y\subseteq A$ a {\em consistent} set of arcs if it contains exactly one path from origin to destination in $G_i$, for each $i\in T$, and exactly one disjunctive arc for each potential conflict associated with the selected paths. 

Let $Y$ be the family of all consistent sets of arcs in $A$, and, for $y\in Y$, let $G(y) = (N(y),y\cup A^Q)$ be the subgraph of $G$ obtained from $G$ by deleting the arcs $(A^R\cup A^C) \setminus \{y\}$ and all isolated nodes. A (feasible) {\em schedule} for $G(y)$, is a vector $t^y\in \reals^{N(y)}_+$ satisfying $t^y_v - t^y_u \geq l_{uv}$ for every $(u,v)\in A(y)$. 

A directed {\em cycle} $K$ of $G(y)$ is called a {\em positive} cycle if $l(K) = \sum_{(u,v)\in C} l_{uv} > 0$.   The following is a very well known result in network flow theory (see, e.g., \cite{ahujia1993network}): 
\begin{lemma}\label{le-acyclic}
A (feasible) schedule $t^y$ for $G(y)$ exists if and only if $G(y)$ does not contain 
positive cycles.  
\end{lemma}

If $G(y)$ does not contain a positive cycle, then it admits a longest path tree $H(y)$ from $s$ to all other nodes. Let $l^y_u$ be the length of a longest path from $o$ to $u$ in $G(y)$, for $u\in N(y)$. One can show that $t^y_u = l^y_u$, for $u\in N(y)$, is a feasible schedule for $G(y)$, and we denote such schedule as $t(y)$. 

We now have an ideal procedure to find a feasible routing and a schedule  for the trains $i\in T$. 
\begin{itemize}
    \item Find a consistent set of arcs $y\in Y$ (corresponding to choosing a route for each train and a winner for each conflict associated with such routing), such that $G(y)$ does not contain positive cycles. 
    \item Compute a longest path tree $H(y)$ in $G(y)$ and let $t^y_u =$ the length of a longest path from $o$ to $u$ in $G(y)$, for $u\in N(y)$. 
\end{itemize}

In other words, if $G(y)$ does not contain positive cycles, the consistent set of  routing and conflict arcs $y$ can be associated with the unique schedule $t(y)$, and so we immediately have a plan. We  now look at the quality of such plan. To this end, we need to introduce the cost $c(t)$ of a schedule $t$. In train dispatching problems like UP's, we are given target arrival times of the trains at destination(s), and the cost of a schedule increases with the delays of the trains with respect to the target times. So, the cost function $c$ is monotone non-decreasing with $t$. With such objective function, one can easily show that, given a vector $y\in Y$, and assuming $G(y)$ does not contain a positive cycle, then the schedule $t(y)$ associated with a longest path tree $H(y)$ in $G(y)$ is actually a minimum cost schedule for the train dispatching problem.  
 Then the dispatching problem can be restated as

\begin{problem}\label{prob:dispatching} {\em\bf [Dispatching Problem]}
Given a dispatching graph $G=(N,A^R\cup A^C \cup A^Q)$, with arc lengths $l\in\reals^A$ and cost function $c: \reals^N \rightarrow \reals$ monotone non increasing, find a consistent $y\in Y\subseteq A^R\cup A^C$ such that $G(y)$ does not contain a positive cycle and $c(t(y))$ is minimized.
\end{problem}

On the other hand, if for each  $y\in Y$ we have that $G(y)$ does contain (at least) a positive cycle, then the current instance of Problem \ref{prob:dispatching} is infeasible, and we say that the trains of $T$ are {\em bound-to-deadlock}. In this case, there is no way we can route all the trains to their final destinations as every combination of feasible trajectories for the involved trains ends up in a deadlock \cite{ArbibEtAl90}. Consequently, we can define the following 
\begin{problem}\label{prob:bound-to-deadlock} {\em\bf [Bound-to-Deadlock Problem]}
Given a dispatching graph $G=(N,A^R\cup A^C \cup A^Q)$, with arc lengths $l\in\reals^A$, check whether there exists a consistent $y\in Y\subseteq A^R\cup A^C$ such that $G(y)$ does not contain a positive cycle.
\end{problem}

Clearly, the above two problems could be merged in a single definition, but we prefer to keep them distinct as they give rise to different functionalities of the AMP. 

\section{Features of the system} 
 
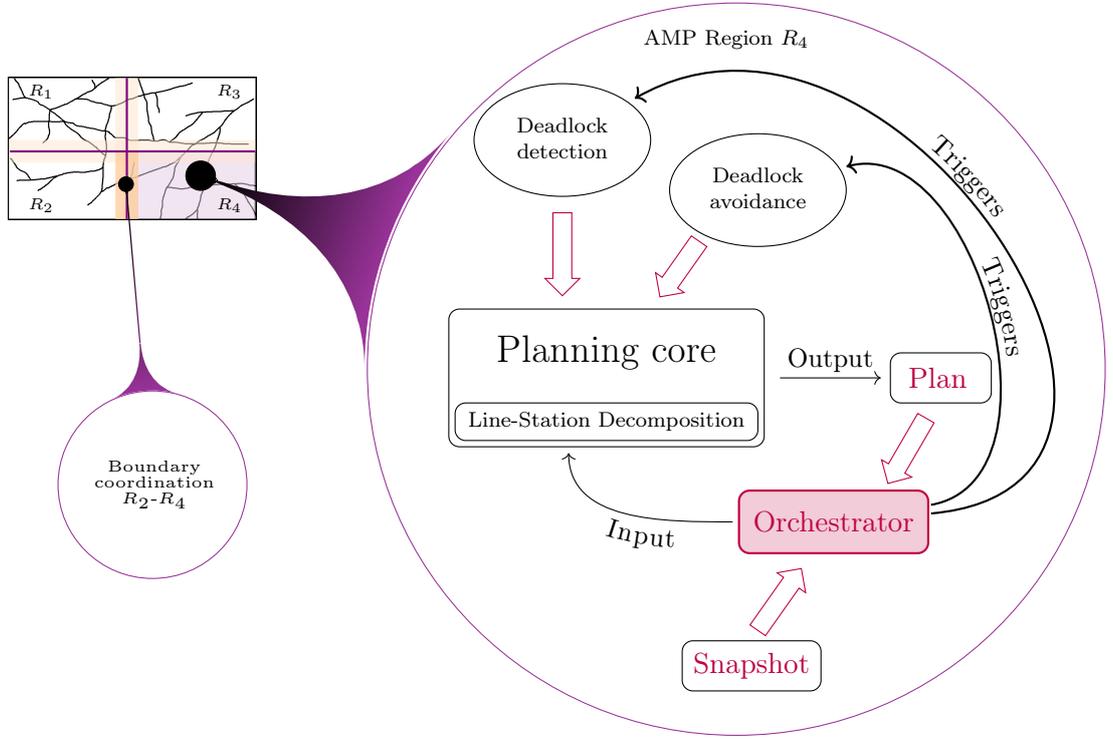
\begin{figure}
    \centering
    \begin{tikzpicture}[scale=0.75]
    \node[](rail)at (0.5,1.1){\includegraphics[scale=0.15]{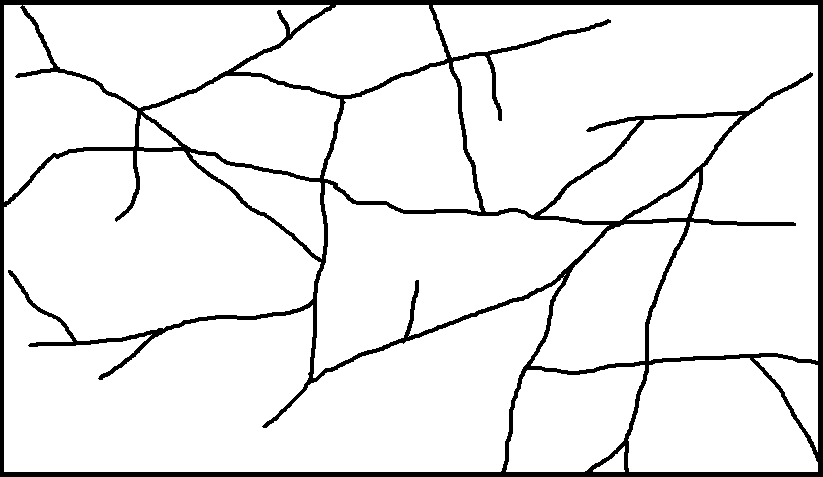}};
    %fill 1 region
    \draw[fill=violet!20, opacity=.5] (0.4,-0.15)--(0.4,1.05)--(2.65,1.05)--(2.65,-0.15)--cycle;
    %borders
    \draw[line width=3mm, orange!20, opacity=.5] (0.4,2.35)--(0.4,-0.15);
    \draw[line width=3mm, orange!20, opacity=.5] (-1.65,1.05)--(2.65,1.05);
    \draw[line width=3mm, orange!50, opacity=.5] (0.4,1.055)--(0.4,-0.15);
    \draw[violet, thick](0.4,2.35)--(0.4,-0.15);
    \draw[violet, thick](-1.65,1.05)--(2.65,1.05);
    %regioni
    \node[](R1) at (-1.1,2.1){\fontsize{6pt}{10}\selectfont {\bf $R_1$}};
    \node[](R2) at (-1.1,0.1){\fontsize{6pt}{10}\selectfont {\bf $R_2$}};
    \node[](R3) at (2.2,2.1){\fontsize{6pt}{10}\selectfont {\bf $R_3$}};
    \node[](R4) at (2.2,0.1){\fontsize{6pt}{10}\selectfont {\bf $R_4$}};
    %mindmap
    \node[]at (9.5,-2.8){\begin{tikzpicture}[rotate=-100, grow cyclic, text=white, level distance = 75mm]
        \node[concept, text width=17mm, fill=black, scale=0.2] {}
	    [clockwise from=80]
        child[concept color=violet!80,text=black]{node(2)[concept, text width=17mm, scale =5]{}};
    \end{tikzpicture}};
    \node[]at (0.85,-2.95){\begin{tikzpicture}[rotate=-165, grow cyclic, text=white, level distance = 40mm]
        \node[concept, text width=17mm, fill=black, scale=0.1] {}
	    [clockwise from=80]
        child[concept color=violet!80,text=black]{node(2)[concept, text width=17mm, scale =1.2]{\fontsize{5pt}{5}\selectfont \begin{tabular}{c}Boundary\\ coordination\\$R_2$-$R_4$\end{tabular}}};
    \end{tikzpicture}};
    %struttura regionale
    \node[]at (13, -1.9){\begin{tikzpicture}[scale=0.83]
    \node[] (AMP) at (3.8,5) {\fontsize{8pt}{10}\selectfont AMP Region $R_4$};
     \draw[rounded corners](-0.6,-1.5)--(4.4,-1.5)--(4.4,0.7)--(-0.6,0.7)--cycle;
     \node[] (PCORE) at (1.9,0) {\fontsize{14pt}{10}\selectfont Planning core};
     \draw[rounded corners](3.1,-5.4)--(3.1,-4.6)--(5.3,-4.6)--(5.3,-5.4)--cycle;
     \node[purple] (D) at (4.2,-5) {\fontsize{11pt}{10}\selectfont Snapshot};
     % \draw[rounded corners, draw=purple!15, fill=purple!5](0.02,0.68)--(1.6,0.68)--(1.6,0)--(0.02,0)--cycle;
     %\node[] (D) at (0.8,0.35) {\fontsize{5.5pt}{10}\selectfont Clock: 10sec};
     \draw[rounded corners](-0.5,-1.4)--(4.3,-1.4)--(4.3,-0.8)--(-0.5,-0.8)--cycle;
     \node[] (Dec) at (1.9,-1.1) {\fontsize{8pt}{10}\selectfont Line-Station Decomposition};
     %\draw[rounded corners](1.5,-2.9)--(4.5,-2.9)--(4.5,-3.9)--(1.5,-3.9)--cycle;
     \draw (1.2,3.4) ellipse (1.4cm and 0.9cm);
     \node[] (OVER) at (1.2,3.4) {\fontsize{8pt}{10}\selectfont \begin{tabular}{c}Deadlock \\detection \end{tabular}};
     \draw (4.3,2.6) ellipse (1.4cm and 0.9cm);
     \node[] (DIS) at (4.3,2.6) {\fontsize{8pt}{10}\selectfont \begin{tabular}{c}Deadlock \\avoidance \end{tabular}};
     \draw[rounded corners](6.4,0)--(8,0)--(8,-0.8)--(6.4,-0.8)--cycle;
     \node[purple] (PLAN) at (7.15,-0.4) {\fontsize{11pt}{10}\selectfont Plan};
     \draw[rounded corners, thick, purple, fill=purple!20](4,-3.2)--(4,-2.2)--(7,-2.2)--(7,-3.2)--cycle;
     \node[purple] (P) at (5.5,-2.7) {\fontsize{11pt}{10}\selectfont Orchestrator};
     %frecce
     \node[single arrow, draw=purple, fill=white, minimum width = 10pt, single arrow head extend=3pt, minimum height=10mm, rotate=55] at (4.6,-4) {};
     %\draw[->, thick, purple, shorten <=0.1cm, shorten >=0.2cm] (D) to (P);
     \node[single arrow, draw=purple, fill=white, minimum width = 10pt, single arrow head extend=3pt, minimum height=10mm, rotate=-120] at (6.7,-1.5) {};
     %\draw[->, thick, purple, shorten <=0.3cm, shorten >=0.25cm] (OVER) to (PCORE);
     \node[single arrow, draw=purple, fill=white, minimum width = 10pt, single arrow head extend=3pt, minimum height=11mm, rotate=-90] at (1.2,1.65) {};
     %\draw[->, thick, purple, shorten <=0.2cm, shorten >=0.4cm] (DIS) to (PCORE);
     \node[single arrow, draw=purple, fill=white, minimum width = 10pt, single arrow head extend=3pt, minimum height=9mm, rotate=-125] at (3.1,1.4) {};
     %\draw[->, thick, purple, shorten <=0.1cm, shorten >=0.25cm] (PLAN) to (P);
    \draw[<-, thick, shorten <=0.2cm, shorten >=0.1cm, postaction={decorate,decoration={raise=4.1pt, text along path,text align=center, text={Triggers{\color{white}{}}}}}] (DIS) to[out=15, in=10, looseness=1.1] (P);
     \draw[<-, thick, shorten <=0.15cm, shorten >=0.1cm, postaction={decorate,decoration={raise=4.1pt, text along path,text align=center, text={Triggers{\color{white}{AAA}}}}}] (OVER) to[out=30, in=5, looseness=1.6] (P);
     \draw[<-, postaction={decorate,decoration={raise=-10pt, text along path,text align=center, text={{\color{white}{...}}Input}}}]  (1.3,-1.6) to[out=-90, in=180, looseness=1] (3.9,-2.7);
     \draw[->, postaction={decorate,decoration={raise=4.1pt, text along path,text align=center, text={{\color{white}{}}Output}}}]  (4.65,-0.4) to (6.25,-0.4);
    
    \end{tikzpicture}};
    \end{tikzpicture}
    \caption{An overview of the system and its continuous planning process.}
    \label{fig:AMP}
\end{figure}

An overview of the system with its features and planning process is given in Figure \ref{fig:AMP}. 
 
A major challenge in train dispatching, like for any real-time application, is precisely the dynamic nature of the problem. To solve the problem effectively one should ideally have a representation of the input that adheres as much as possible to reality. For train dispatching, this primarily means accurately representing the overall state of the network, such as train positions, signals, interlocking, blocks, etc, at a given moment in time. We'll refer to the particular configuration in time and space as a {\em snapshot} of the network. Now, provided such input, like for any algorithm, some computational effort is necessarily required to produce an output. We'll refer to this as a {\em planning cycle}, namely representing all the computational activity that occurs between the reference moment in time (associated with the snapshot) and the moment in which the cycle is declared over (which should coincide with the successful production of a plan). The challenge, again, stems from the real-time nature of the problem. The greater the length of the planning cycle, the higher the risk of producing a suboptimal or indeed effectively infeasible plan, i.e. a plan that is not compatible with the {\em current} state of the network, since during the time elapsed the state of the network may in general have changed from the snapshot. On the other hand, a longer planning cycle allows more time for an algorithm to find good, possibly optimal solutions to the problem. Thus, the trade-off.\footnote{The planning cycle duration in AMP is set to 10 seconds by default, however depending on the application can be configured accordingly}

Now we give an overview of how the different components of the system interact during a given planning cycle. The complexity of the real-time application combined with the number of different aspects of the problem that have to be dealt with requires a dedicated "orchestration" layer that receives the ever-changing information from the field, schedules the different components and ensures the right flow of data between them.  
This layer is in charge of taking the snapshot of the network that kicks off the start of a planning cycle. At this point, the Planning core ({\em Planner}) is executed with the current state of the network and the output (from the previous planning cycle) of the various ancillary components dedicated to deadlock detection, deadlock avoidance and more (described later in the section). The planning cycle ends once the timeout is reached, and the best plan produced is published to the TMS, where it automatically translates into signal and interlocking commands that are sent to the field.\footnote{By default the system functions in fully automatic mode, however on specific areas or trains it can also be operated in "decision support" mode, where the commands must be manually confirmed by the dispatcher} Such plan however also becomes an input for the deadlock detection and avoidance modules at the {\em following} planning cycle. Indeed, these are triggered based on specific logic that compares the plan with new events that have occurred in the field during the planning cycle (and therefore were not accounted for as part of the snapshot). If these events are deemed definitely or even potentially incompatible with the published plan, the relevant modules are executed and, similarly to the {\em Planner}, use the remainder of the planning cycle to produce a solution to their specific problem. This output becomes one of the inputs of the {\em Planner} at the following planning cycle and so forth. 

This very briefly summarizes the planning process at the heart of AMP, which carries on continuously, in a loop, while the system is up and running. 
Figure \ref{fig:AMP} additionally references how AMP can be run on multiple coordinated regions, a topic discussed later in a section dedicated to scalability. In this example the controlled area is divided into 4 planning regions ($R_1, \dots R_4$), where each region is managed by an instance of AMP, and a boundary coordination algorithm is applied to each border to ensure the smooth handover of trains. An integer programming model to construct suitable partitions of the UP network into planning regions is described in \cite{Ve24}. 

\subsection{Train dispatching by Branch\&Bound}\label{sec:planner} 
This functionality is handled by the core component of the system, the {\em Planner}. The goal of this component is to find a solution to a specific instance $G = (N,A)$, $l\in \reals^A$, $c:\reals^N \rightarrow \reals$ of train dispatching, namely return a consistent set of arcs $y\in Y \subseteq A^R\cup A^C$, with associated longest path tree schedule $t(y)$. This is carried out by means of a specialized truncated branch and bound method (B\&B), called {\em continuous rolling horizon} (CRH). 
The CRH scheme is inspired by a classical heuristic decomposition approach, 
called Rolling Horizon. In this framework, the planning horizon H is subdivided (often partitioned) 
into intervals $H_1, H_2, \dots, H_q$. Then a sequence of subproblems associated with a reduced planning 
horizon $H_1 \cup H_2 \cup \dots \cup H_i$ is solved, for $i \leq q$. When solving the problem relative to $H_1$, we consider 
only the vector $x_1$ of decision variables associated with events which are occurring before the end of $H_1$. All other decision variables are removed from the current problem, which is then solved 
(possibly to optimality). Let $\xi_i$ be the solution found. At the next iteration, the subproblem 
associated to the planning horizon $H_1 \cup H_2$ is built, with decision variable vectors $x_1, x_2$. However, 
we let $x_1 = \xi_1$, namely the decision variables associated with the previous subproblem are fixed (all 
of them or just a subset in some versions) to their previous best value. The method iterates and, 
at iteration $i$, we solve the problem corresponding to $H_1 \cup H_2 \cup \dots \cup H_i$, with (as subset of) the 
variables fixed to some precalculated values. In vehicle scheduling this approach is followed in 
many articles (e.g. \cite{Da09, Ni12, Sa13, Zh16} and many others). The idea is certainly appealing, 
especially in cases when predictions on vehicles behaviour become fuzzier as time passes. 
However, its actual implementations suffer from various downsides. First, it may happen that, at 
a certain iteration, the current subproblem is infeasible. Then, unless we have a suitable 
backtracking mechanism, we are not able to decide whether the overall problem is infeasible, or 
the deadlock is caused by some previous, incautious fixing of decision variables. Most rolling 
horizon algorithms do not include any backtracking mechanism. In other cases, the authors are 
vague on the question, or they invoke an external action - as in \cite{Sa13}, where human operators are asked to resolve such situation. In addition, no mechanism is suggested to improve the 
solutions, to compute bounds, etc. Similarly, no ways to produce lower bounds are suggested and 
implemented. Another difficult matter is the determination of intervals $H_1$, $H_2$, $\dots$, $H_q$. If the 
intervals are few, we risk tackling too large instances. On the other hand, when they are too many, 
the solutions can be poor and the risk of generating infeasible subproblems is very high. By 
combining the rolling horizon paradigm with a B\&B scheme, the CRH approach is able 
to overcome all of the above-mentioned limitations. The method is exact and can both prove the 
optimality of the solution produced or prove that the instance is infeasible. 

\textcolor{black}{
We  briefly sketch how our B\&B proceeds. Any B\&B method (see \cite{nemhouserWolsey88}) requires a few main ingredients: an enumeration scheme to recursively build solutions, a (lower) bounding procedure, an incumbent solution $y^I$ whose cost $c(t(y^I))$ defines an upper bound $UB$ on the value of the optimal solution to the problem. In our case, a feasible solution is a consistent selection of arcs $y\in Y\subseteq  A^R \cup A^C = A^0$. We (implicitly) build all possible selections by starting with  an empty subset $x^0$ of arcs and recursively perform a branching operation at each decision of including/excluding conflict/routing arcs in the current partial solution. Assume that at branching node $k$ we add arc $(u^k,v^k)\in A^{k-1}$, and let $x^k = x^{k-1}\cup \{(u^k,v^k)\}$, $A^k = A^{k-1} \setminus \{(u^k,v^k)\}$. The set $x^k$ is called {\em partial selection}, and we call {\em extension} of $x^k$ any consistent selection $\bar y\in Y$ such that $x^k\subseteq \bar y$. Let $N^k \subseteq N$ be the set of nodes which are reachable from $s$ in $G(x^k)$. If $G(x^k)$ contains a positive dicycle, then we can discard this branching because $x^k$ cannot be extended to a feasible consistent selection. Otherwise, we compute the longest path tree from $s$ to every node in $N^k$, and we let $t^k\in \reals^{N^k}$ be the associated schedule and $c^k$ be the projection of $c$ into $\reals^{N^k}$ ($c$ is separable). If $c^k(t^k) \geq UB$ then we can discard the current branching because $x^k$ cannot be extended to a consistent selection which improves the incumbent. Indeed, since $c$ is monotone non-decreasing, $c(t(\bar y)) \geq c^k(t^k)$ for any extension $\bar y$ of $x^k$.\\ \\
Otherwise, $c^k(t^k) <UB$, and if  $x^k\in Y$, that is $x^k$ it is a consistent selection (and thus  $N^k = N$), then we let $UB = c(t^k)$ and update the incumbent by letting $y^I = x^k$.  
Otherwise we iterate with the next branching. \\ \\
%
%The way the arc $(u_k,v_k)$ is selected at each iteration is what defined %our {\em continuous rolling horizon approach}. 
Note that $G(x^k)$ is obtained from $G(x^{k-1})$ and, at branching node $k$, we have at hand the parent solution $t^{k-1}$. This can be exploited when computing the longest path for $G(x^k)$ in very effective ways. 
In particular, we have developed some specialization of the well-studied incremental techniques for shortest/longest path tree computations on acyclic digraphs (see \cite{bender2015new, haeupler2008faster}). %
Besides the CRH approach, we also make use of a Benders' like decomposition scheme based on decomposing the network into lines and stations, as described in \cite{lamorgese2015exact}. }

\subsection{Deadlock avoidance}
Deadlock avoidance is arguably the first responsibility of a dispatcher. In principle the routing and scheduling solutions found by the {\em Planner} are deadlock-free, that is, should not lead to a deadlock when carried out. However there are special circumstances in which, without further attention, deadlocks may indeed occur. The first is related to network disruptions, i.e. when parts of the network become unavailable for inbound trains. Deadlocks are more prone to occur during disruptions because the capacity of the network is often abruptly reduced and dispatchers face unfamiliar and critical circumstances. AMP’s routing and scheduling core algorithm seamlessly factors in any track blocks or slowdowns in its planning process and accordingly produces optimized plans. However more has to be done under more extensive disruptions. In particular if there is a subset of trains that cannot route around the disruption, these should be “safely” parked somewhere until the disruption is over, i.e. at a location that allows traffic to route around it. The purpose is to avoid trains driving past the last location where they can be routed around, which could otherwise lead to further disruptions and deadlocks.
The second key deadlock risk is related to the time horizon of the routing and scheduling plan. While this time horizon should be chosen suitably large, there is always a theoretical possibility that decisions taken at a given time could lead to deadlocks {\em beyond} the planning horizon. Again, this is particularly the case when many dispatched trains are overlength. To avoid deadlock-inducing decisions for this subset of critical trains, it can be beneficial to take a look beyond the planning horizon. 
%Because of this, it is crucial to ensure that the choices made by the planner will not lead to a deadlock beyond the planning horizon (notice that the planner itself has the ability to avoid deadlocks within its planning horizon). To prevent deadlocks further ahead, there are two main aspects to consider: firstly, that when disruptions happen on the railway network, trains cannot be planned up to the disrupted location without taking into account the presence of other trains piling up on the other side of the disruption. Hence, these trains need to be parked in such a way that, when the disruption is lifted, planning can resume.
%These parking locations are called \emph{safeplaces}. Secondly, that whenever longer trains, referred to as \emph{overlength trains}, are circulating, their only meet location may be beyond the planning horizon of one of the trains. Hence, it is possible that the planner pushes the other trains further ahead and creates a deadlock. 
These issues are tackled in AMP respectively with a module that assigns \emph{safeplaces} and a module that deals with overlength trains.
    
\subsubsection{Safeplacer module} 
Safeplaces are defined as locations where a train can be stopped without impeding the transit of other trains. Notice that safeplaces are train-dependent as, due to the different length of the trains, a location that could be used as safeplace for a train may be not a valid location for another (longer) train. As explained before, safeplaces are only needed in hazardous situations, where extra care is necessary to prevent and manage deadlocks. In particular, a call to the Safeplacer is triggered when some trains are reaching either a disruption on the network or a subset of deadlocked trains. Trains in need of a safeplace are grouped based on the triggering event, then for each group of trains this module assigns them safeplaces so to guarantee that the {\em SP-protocol} defined by UP is satisfied. 

The safe place assignment (with no decision support) is usually decided by experienced dispatchers in a "greedy" fashion. 
The automatic procedure embedded in AMP is based instead on a ILP approach to the problem that constructs optimal solutions in fractions of seconds, providing strong support to dispatchers in their task.
In finding the best possible holding locations, i.e. the safeplaces, the SP-protocol pursues three distinct (and somehow conflicting) objectives:
O1) The primary goal is to avoid that when the disruption is lifted trains are bound-to-deadlock; O2) Secondly, to leave some paths free of trains so that worker trains can reach the disrupted area; O3) Finally, to "push" trains as far along their routes as possible, so as to mitigate the impact of the disruption in terms of train delays.

Objectives O1 and O2 are implemented as constraints of the ILP formulation. 
In contrast O3 is expressed in the objective function of the formulation. 
As reported in detail in \cite{croella2022disruption}, the overall formulation is mostly defined by vertex-packing type constraints, plus some extra constraints needed to "activate" the alternative paths that model objective O2.  
The results of the computational experiments that prove the effectiveness of such formulation on a data set of realistic instances are also presented in  \cite{croella2022disruption}. 
In particular, observe that almost all the instances could be solved to optimality in a few hundredths of a second and none of them requires more than 0.6 seconds of CPU time.

\subsubsection{Overlength train module}
AMP also includes a module dedicated exclusively to dispatching “overlength” trains, which tend to be common in freight traffic and, as mentioned, present a major deadlock hazard. The need presents itself specifically when the  {\em Planner}’s “regular” planning horizon (generally set to a few hours) is insufficient to identify a deadlock risk coming from overlength trains that are instead still far apart. To tackle this issue, the overlength train module "looks" beyond the regular planning horizon and helps the {\em Planner} to avoid deadlock-inducing decisions for this subset of critical trains. 

The overlength trains module is launched based on two criteria: the first is the maximum time that should elapse between two consecutive runs; the second is the occurrence of events that could potentially make previous decisions not feasible or desirable anymore. Interestingly, the algorithm underlying this module is the CRH itself, simply executed with a broader time horizon. The consequently increased computational complexity is tackled by focusing only on a subset of suitably chosen trains. These amount to (a) all overlength trains (typically defined as trains longer than a certain critical length) and (b) a subset of the remaining trains that are positioned at locations that make them a potential hindrance for smoothly dispatching the overlength trains. 
Once the subset of trains has been determined, an instance of the CRH with a suitably large time horizon is executed. Relevant scheduling and routing decisions for the overlength trains are extracted from the resulting plan and passed on to the Planner as fixed decisions. 
This decomposition approach has proven successful to strongly mitigate the risk of deadlocks induced by overlength trains.

\subsection{Deadlock detection}

Deadlocks are usually the result of unexpected disruptions combined with a shortsighted view of the effect of dispatching decisions on traffic down the line. The chance of creating deadlocks is significantly higher when many circulating trains are overlength and the locations for routing one past the other are few. Indeed, it happens from time to time that a deadlock occurs but is difficult to spot for dispatchers due to trains being still far from each other. Needless to say, it is of crucial importance to quickly detect deadlocks, as it may still be possible to revert deadlock-inducing commands or, at least, avoid trains travelling too far forward (as solving a deadlock typically requires pulling back one or more trains). Furthermore, in these cases the {\em Planner} component will not be able to find a feasible plan. Indeed, we already showed how an instance $G=(N,A)$ of Problem \ref{prob:dispatching} is infeasible if and only if, for any consistent set of arcs $y\in Y$, $G(y)$ contains at least a positive directed cycle. 
In principle, we could solve Problem \ref{prob:bound-to-deadlock} and prove this condition by running the Planner for sufficient time to explore the full branching tree. However, this is often impossible even for medium-size instances given the harsh computing-time limitations and the well known "curse of exponentiality". Hence we developed ad-hoc algorithms to tackle this issue and attempt to promptly identify the potential for deadlocks.

The procedure we developed is based on two assumptions:
\begin{enumerate}[label=A\arabic*]
\item \label{DDass1} the plan returned by the Planner cannot "generate" a deadlock beyond the time horizon,
\item \label{DDass2} the output of the previous planning cycle, starting at $t_{i-1}$, is a valid plan with routing choices and schedule for all trains up to $t_{i-1} + TH$.
\end{enumerate}

Assumption~\ref{DDass1} is granted by the deadlock avoidance modules embedded in the system's components, as described in the previous section. Assumption~\ref{DDass2} holds when the whole system is in operation, as the Planner publishes a valid plan every few seconds.
Based on these assumptions, at time $t_i > t_{i-1}$ a deadlock may occur only if a potentially critical event that disrupts the previous plan happened since the last Planner cycle was performed, i.e. between time $t_{i-1}$ and $t_i$. Events of this type, such as, e.g.,  a track blockage, a new train entering the network or simply a train extended with some cars, are called the \emph{triggers}. When any such triggers occur, before running a new planning cycle, we perform a feasibility check by invoking the deadlock detector.  

First, we identify a subset of trains that is most likely to be affected by the trigger and the corresponding area. %(the {\em bubble}). 
This identification is carried out heuristically, starting from the train(s) or location(s) "nearby" the location where the trigger occurred. 

% The starting point is, of course, the train that underwent the critical change of schedule. The main idea we adopted is the following: if we are still able to route this train to a location where it does not impede the traffic of the other trains, then no train is bound to deadlock. Hence we initially consider the portion of the network between the current train location and its first safeplace. The trains that are located inside this area are added to the subset of potentially critical ones. Afterwards, this subset of trains may be broadened if the station selected as safeplace location was also the meet or pass designed location for other pairs of trains. In order to keep the size of the instance as small as possible, these trains are not added to the subset if we can easily find them a safeplace outside of the initial area.

Then, we apply an ad-hoc algorithm in order to identify whether a subset of these trains is bound-to-deadlock or not. When such trains are detected, they are notified to the Planner, that holds them at their current locations and keeps planning for the remaining trains. %the bubble contains bound-to-deadlock trains. 
%\CAR{Questo deve essere scritto un po' meglio: chi sono i treni coinvolti? Vengono tutti fermati nella loro locazione attuale? Oppure solo alcuni?: 

%... there are deadlocked trains among the selected ones or not. When a deadlock is detected, the involved trains are notified to the Planner, that holds the deadlocked trains at their current locations and keeps planning for the remaining ones ...} 

In the remaining of this section we sketch the two different deadlock-detectors modules embedded in AMP. 

\subsubsection{Two-train deadlock detector \cite{2trains}} 
Deadlocks involving only a pair of trains are, perhaps surprisingly, the most frequent in freight railroad (see \cite{Pachl2011}). A dedicated module in AMP implements the algorithm developed in \cite{2trains} in order to check whether two trains, each with a given starting position and final destination, are bound-to-deadlock. We refer to this as the {\em bound-to-deadlock 2-train problem}. 
Using a combinatorial result by McCormick et al. \cite{LiMcCormick},  the following has been proven in \cite{2trains}: 
\begin{theorem}
The Bound-to-Deadlock 2-Train Problem can be solved in polynomial time for any network and any train length.
\end{theorem}
Moreover, in \cite{2trains} we also develop a pseudo-polynomial algorithm that can solve very efficiently all our real-life instances.

This module is invoked at the beginning of each planning cycle to detect whether any events have created a two train deadlock. In practice it is not necessary to repeat the feasibility check for every pair of trains. It suffices to focus only on a subset of these, namely those affected by a deadlock trigger or trains that were too far apart to be covered by the {\em Planner } within its time horizon.

\subsubsection{Multiple-train deadlock detector}
We solve  Problem \ref{prob:bound-to-deadlock} for an arbitrary number of trains by means of a binary integer program which extends the feasibility constraints of the Path\&Cycle formulation introduced in \cite{lamorgese2019noncompact}. For sake of brevity, we present here only the most important features of the model. 

We are given $T$ the set of trains, the movement graph $G_i = (N_i, A_i)$, with $o_i, d_i \in V_i$, for each $i \in T$, and the dispatching graph $G = (N, A^R\cup A^C \cup A^Q)$. 
For each $i \in T$ and $a \in A_i$, we use the binary variable $z^i_a$ that is 1 if train $i$ uses the (block section associated with) arc $a$ (i.e. $\bf z^i$ is the incidence vector of the route chosen for $i$). Next, we need to model the fact that block sections are contended by trains. This imposes that certain pairs of movements cannot happen simultaneously. Such pairs are represented precisely by pairs of arcs in the disjunctive arcs set $A^C$. In particular, for a pair of distinct trains $i,j\in T$ and some pairs of routing arcs $a_i = (u_i,v_i) \in A_i$ and $a_j= (u_j,v_j)\in A_j$, either $i$ runs through the block section corresponding to $a_i$ before $j$ runs through the block section corresponding to $a_j$, or viceversa. This implies that $(v_i,u_j)$ and $(v_j,u_i)$ belongs to $A^C$. We call the pair $\{a_i,a_j\}$ an {\em incompatible couple} and we let $J$ be the set of all such couples.
For each $\{a_i, a_j\}\in J$ (associated with the pair of trains $i,j$), we introduce the binary variables $x_{a_i, a_j}$ and $x_{a_j,a_i}$ and $x_{a_i,a_j}$ ($x_{a_j,a_i}$) assumes value 1 if $i$ traverses $a_i$ before $j$ enters $a_j$ ($j$ traverses $a_j$ before $i$ enters $a_i$, resp.), and 0 otherwise. As already explained, this amounts to selecting one of the two arcs\footnote{Remark that  $a_i,a_j\in A^R$ whereas $(v_i,u_j),(v_j,u_i)\in A^C$} $(v_i,u_j), (v_j,u_i)\in A^C$  in the disjunctive graph. 

% for a pair of distinct trains $i,j\in T$ and some pairs of arcs $a_i \in A_i$ and $a_j\in A_j$, either $i$ runs through the block section corresponding to $a_i$ before $j$ runs through the block section corresponding to $a_j$, or viceversa. We call the pair $\{a_i,a_j\}$ an {\em incompatible couple}, and we denote by $J$ the set of all incompatible couples (for all trains). Then, for each incompatible couple $\{a_i, a_j\}\in J$ (associated with the pair of trains $i,j$), we introduce the binary variables $x_{a_i, a_j}$ and $x_{a_j,a_i}$ and $x_{a_i,a_j}$ ($x_{a_j,a_i}$) assumes value 1 if $i$ uses $a_i$ before $j$ enters $a_j$ ($j$ uses $a_j$ before $i$ enters $a_i$, resp.), and 0 otherwise.

Next we present the (linear) constraints of the formulation, starting from the set of inequalities that ensure that $\bf z$ is, for each train, the incidence vector of a origin-destination path in the associated movement graph. These are the standard {\em node balance} constraints for flow problems (see \cite{ahujia1993network}):
\begin{eqnarray}
    \sum_{a \in \delta^+(u)} z^i_{a} = \sum_{a \in \delta^-(u)} z^i_a &  ,\forall \, u \in N_i \setminus \{o_i, d_i\}, \, \forall \, i \in T, \label{eq:path1}\\
    \sum_{a \in \delta^+(o_i)} z^i_a = 1, &  \forall \, i \in T,\label{eq:path2}\\
    \sum_{a \in \delta^-(d_i)} z^i_a = 1, &  \forall \, i \in T.\label{eq:path3}
\end{eqnarray}

Next, we have the incompatibility constraints, preventing two trains to use contended block sections at the same time. Therefore, if  $\{a_i, a_j\}$ belongs to $J$, the path of train $i$ uses $a_i$ and the path of train $j$ uses $a_j$, then either $t_i$ uses $a_i$ before $t_j$ uses $a_j$ or viceversa. Then, exactly one variable between $x_{a_i, a_j}$ and $x_{a_j, a_i}$ must take value $1$. Otherwise both variables are 0. This condition is enforced by the following constraints: 
\begin{eqnarray} \label{eq:yactivate01}
    x_{a_i,a_j} + x_{a_j,a_i} & \geq & z^i_{a_i} + z^j_{a_j} - 1, \qquad \forall \, (a_i, a_j) \in J \label{eq:yactivate1}\\
    x_{a_i,a_j} + x_{a_j,a_i}& \leq & z^i_{a_i}, \qquad \qquad \qquad \;  \forall  \, (a_i, a_j) \in J \label{eq:yactivate2}\\
    x_{a_i,a_j} + x_{a_j,a_i}& \leq &z^j_{a_j}, \qquad \qquad \qquad \; \forall   \, (a_i, a_j) \in J
    \label{eq:yactivate3}
\end{eqnarray}

\begin{comment}
\CAR{Questo  vincolo è del tutto incomprensibile. La variable x con un indice (su un nodo!) non esiste, halted non si sa che vuol dire qua, eccetera. Meglio forse rimuoverlo: 

Moreover, each train $t$ precedes all the others in its initial position (since it occupies it at the current time). This means that, if $t$ is halted,  in $o_t$, all the other trains cannot use the resources blocked by $o_t$.
\begin{eqnarray}
     z_{a_i} + x_{o_t} \leq  1 &\forall \, t \in T, \forall \, (a_i, a_j) \in J |\, a_j \in \delta^+({o_t}), \label{eq:y_iniz1}
\end{eqnarray}
Instead, if $t$ moves forward, all other trains must use the resources blocked by $o_t$ only after $t$.
\begin{eqnarray}
     y_{a_i, a_j} & =  0 &\forall \, t \in T, \forall \, (a_i, a_j) \in J |\, a_j \in \delta^+(o_t). \label{eq:y_iniz2}
\end{eqnarray}

}
\end{comment}
Observe now that  vectors  $\bf z$ and  $\bf x$ form the incidence vector of a selection $y \subset A$ of arcs of the dispatching graph $G$. We know by Lemma \ref{le-acyclic} that, to be feasible (i.e., trains not bound-to-deadlock), the subgraph $G(y)$ of $G$ associated with such selection must not contain positive directed cycles. Then, denoting by $\mathcal K_+$ the set of positive directed cycles of $G$, the following constraint must be satisfied for every $K\in \mathcal K_+$: 
\begin{eqnarray} \label{eq:cycle}
    \sum_{a_i \in K\cap A^R} z^i_{a_i} + \sum_{({v}_i, {u}_j) \in K\cap A^C}  x_{a_i, a_j} \leq |K| - 1 
\end{eqnarray}
We solve the overall 0,1-linear formulation by iterative delayed row generation.
In particular, since the number of cycles in a graph grows exponentially with the size of the graph, so does the number of constraints \ref{eq:cycle}.
Therefore, we separate them in a cutting plane fashion during the branch-and-cut procedure used to solve the formulation.

% Observe now that, since the $z$ variables satisfy the constraints \eqref{eq:path1}-\eqref{eq:path3}, they define a set $\cal P$ of origin-destination paths, one for each train in $T$. Therefore, they define a subgraph $G^\star$ of the dispatching graph $G(\cal P)$ that satisfies conditions \emph{(i)} and \emph{(ii)} of Theorem~\ref{th:BDproblem}. Moreover, condition \emph{(iii)} is implemented by constraints \eqref{eq:yactivate1}-\eqref{eq:yactivate3}. Now, in order to satisfy condition \emph{(iv)}, we need to ensure that the arcs of $D(\cal P)$ selected by variables $y$ and $z$ do not support any cycle. This task can be easily achieved by using the following constraints.

% \begin{eqnarray} \label{eq:cycle}
%     \sum_{a \in A(C)} z_{a} + \sum_{(\bar{a}_i, \bar{a}_j) \in \bar{J}(C)}  y_{a_i, a_j} \leq |A(C)| + |\bar{J}(C)| - 1 \qquad \forall\,  C \mbox{ cycle of } D(\cal P)
% \end{eqnarray}

\subsection{Scalability} \label{scalability}
To ensure AMP's scalability when the controlled network grows beyond tractability, we introduce two levels of decomposition. The first is geographical, and consists in subdividing the network into different regions. The second level of decomposition is applied to the single region and, in a master-worker fashion, decomposes the problem into a higher-level (line) problem and a lower-level (station) problem.

\subsubsection{Geographical decomposition}
This geographical decomposition approach mimics the approach taken in practice by dispatchers (described in the introduction), albeit on a larger scale, since each AMP instance can generally cover the territory of up to tens of dispatchers. Naturally this requires a layer of coordination between the different planning instances, to ensure that the handover between regions occurs smoothly and in a way that is globally beneficial.
To the best of our knowledge, the coordination approaches explored in the literature iteratively solve the single region problems and a coordination problem, until a consistent solution is found (see, for example, \cite{corman2014dispatching} and \cite{luan2020decomposition}). Although, when a global feasible solution exists, these approaches are known to converge, a practical application may suffer from the long computational times required to reach such solution. 

The train transit across region boundaries in AMP is supervised by a dedicated \emph{coordination module}, whose responsibility it is to decide where each train should cross (i.e the routing) and establishing the order in which the trains traverse the border. These precedences are set to guarantee a feasible/deadlock-free handover between regions and, at the same time, factor in the impact on the global objective. The exact scheduling times remain in the hands of the individual regional Planners, and subject to the constraints defined by the coordination module. 
%\LEO{Questo paragrafo seguente potrebbe essere un pò controproducente svelarlo}
%Opposed to the other approaches in the literature, our coordination module builds a coordination area for each pair of adjacent regions. This is possible because of two main characteristics of the region partition in use: (i) there are no borders involving more than two regions at the same time and (ii) the borders are sufficiently separated, so that any train traverses at most one border within the planning horizon.

%\Leo{Sono indeciso se valga la pena arrivare a questo livello di dettaglio}
%An instance of CRH runs on each coordination area with all the trains that are expected to enter it within the time horizon and produces a schedule for each train up to the point where the train leaves the coordination area. Then, the precedences on the border tracks are extracted and fed back to the regional Planners.

The partitioning of the UP network into planning regions has been determined by  a dedicated algorithm, which aims to define regions that are balanced in terms of “planning complexity” while forming boundaries that have suitable properties for the coordination process to take place seamlessly. The algorithm is based on the MILP model for spatial decomposition described in \cite{Ve24}.

%The partition of the network into regions can be defined by the rail operator based on an existing or desired region partition. However AMP also embeds a dedicated algorithm that defines regions that are balanced in terms of “planning complexity” while establishing boundaries that have suitable properties for the coordination process to take place as seamlessly as possible. 

\subsubsection{Line-station decomposition}
The second level of decomposition is applied to the single region. As already mentioned in Section~\ref{sec:planner}, following the approach presented in \cite{lamorgese2015exact}, the problem is divided into a master (line) problem and a worker (station) problem. In the master problem, stations are collapsed into a single (capacited) railway resource and the master solution provides arrival and departure times for each train in each station and is conflict- and deadlock-free at the line level. The slave problem further decomposes into several smaller problems, one for each station. As in the classical Benders' decomposition, the problems communicate through suitable feasibility and optimality cuts in a delayed row generation fashion. 

% Each station is then treated as a small rail network itself, and routing and scheduling choices are optimized locally, while the schedule is reported globally. The effect of such decomposition on the Planning component is a reduction of the size of the branching tree, so that the quality of the incumbent solution can be increased. \Leo{Da sistemare un po quest'ultima parte, in ogni caso Carlo direi che vanno citati i primi paper su OR e TS}

\section{Conclusions}
In this paper we presented the story of a groundbreaking project commissioned by Union Pacific, one of the largest transportation companies in the world, which led to the deployment of the currently only optimization-based traffic management system that allows the full automation of train dispatching in real-time (A-TMS). Achieving this required tackling a number of challenging optimization problems, which in turn demanded the development of suitable mathematical approaches based, among other things, on graph theory, combinatorial optimization, exact and heuristic methods for integer linear programming and decomposition theory. 
AMP, the optimization module that powers this A-TMS, provides the different functionalities that are necessary to dispatch a rail network automatically and effectively, including some that are particularly relevant for (predominantly) freight-based rail networks. Among these, routing and scheduling (the core functionality for train dispatching), deadlock detection, deadlock avoidance, disruption management, and more. AMP was first rolled out on a region of the UP network in May 2021, and fully deployed system-wide in November 2021. Since then, every 10 seconds an AMP instance computes deadlock-free, conflict-free, optimized plans for each region of the network. 

\bibliographystyle{plain}
\bibliography{biblio.bib}
\end{document}